\documentclass[10pt]{article}
\usepackage{amssymb, amsmath}

\begin{document}
\title{An analogue of a van der Waerden's theorem\\
and its application to two-distance preserving mappings}
\author{Victor Alexandrov}
\date{ }
\maketitle
\begin{abstract}
The van der Waerden's theorem reads that
an equilateral pentagon in Euclidean 3-space $\Bbb E^3$ 
with all diagonals of the same length is necessarily 
planar and its vertex set coincides with the vertex
set of some convex regular pentagon.
We prove the following many-dimensional 
analogue of this theorem: for $n\geqslant 2$,
every $n$-dimensional cross-polytope in $\Bbb E^{2n-2}$
with all diagonals of the same length and all
edges of the same length necessarily 
lies in $\Bbb E^n$ and hence is a convex regular cross-polytope.
We also apply our theorem to the study of 
two-distance preserving mappings of Euclidean spaces.
\par
\noindent\textit{Mathematics Subject Classification (2010)}: 
52B11; 52B70; 52C25; 51K05.
\par
\noindent\textit{Key words}: Euclidean space, pentagon,
cross-polytope, Cayley-Menger determinant, 
Beckman-Quarles theorem
\end{abstract}

\textbf{1. van der Waerden's theorem and its many-dimensional analogue.}
In 1970 B.L. van der Waerden has shown that an equilateral and
isogonal pentagon in Euclidean 3-space $\Bbb E^3$ is necessarily 
planar and its vertex set coincides with the vertex
set of some convex regular pentagon.
For more details about this theorem we refer to 
\cite{Wa70} and \cite{Bo73} and references given there.
More recent results related to this theorem may be
found in \cite{BK11}, \cite{FS12}, and \cite{OH13}.

The van der Waerden's theorem may be reformulated as follows: 
\textit{an equilateral pentagon in Euclidean 3-space $\Bbb E^3$ 
with all diagonals of the same length is necessarily 
planar and its vertex set coincides with the vertex
set of some convex regular pentagon.}
The aim of this paper is to find a many-dimensional 
analogue of this statement by replacing a pentagon by
a polyhedron which, under some conditions
on the lengths of its edges and diagonals, is `surprisingly flat.'
The latter means that we start with a polyhedron located 
in $\Bbb E^N$ and prove that it necessarily lies in
$\Bbb E^n$ for some $n<N$.

Let $n\geqslant 1$ be an integer and 
$\{ {\bf e}_1, {\bf e}_2, \dots, {\bf e}_n \}$ 
be an orthonormal basis in Euclidean $n$-space $\Bbb E^n$.
Denote by $V_n$ the set of the end-points of the vectors 
$\pm{\bf e}_1, \pm{\bf e}_2, \dots, \pm{\bf e}_n $. 
The convex hull of $V_n$ is called the {\it standard $n$-dimensional 
cross-polytope} in $\Bbb E^n$ and is denoted by $S_n$.
Obviously, $V_n$ is the set of vertices of $S_n$.

Let $N$ and $n$ be positive integers.
In this paper, every injective mapping $f:V_n\to\Bbb E^N$
is called an {\it $n$-dimensional cross-polytope} in $\Bbb E^N$.
We prefer this definiton for the following two reasons:

(i) using $f$, the reader may easily reconstruct 
an abstract simplicial complex $C$ which is combinatorially equivalent to $S_n$ 
and whose vertex set is $f(V_n)$ (it suffice to put by definiton
that a simplex $\Delta$ (of arbitrary dimension) with vertices
$w_1,w_2,\dots,w_k\in f(V_n)$ belongs to $C$ if and only if 
the simplex with vertices 
$f^{-1}(w_1),f^{-1}(w_2),\dots,f^{-1}(w_k)$ is a face of $S_n$)
and 

(ii) from our definition the reader may easily see that
the realization of the abstract simplicial complex $C$ 
described in (i) may be degenerate (e.g., the vertices
$w_1,w_2,\dots,w_k\in f(V_n)$ may lay on a single line in $\Bbb E^N$).

For every two points $u,v\in V_n\subset S_n$ there are only 
two possibilities: either $u$ and $v$ are joint together 
by an edge of $S_n$ or $u+v=0$. In the first case the
straight-line segment with the points $f(u)$ and $f(v)$ in
$\Bbb E^N$ is called an {\it edge} of the $n$-dimensional 
cross-polytope $f:V_n\to\Bbb E^N$, while in the second case 
this segment is called a {\it diagonal} of $f$.

The main result of this paper is the following theorem.

\textbf{Theorem 1.} \textit{Let $n\geqslant 2$ be an integer,
let $a$, $b$ be positive numbers, and let $f:V_n\to \Bbb E^{2n-2}$
be an $n$-dimensional cross-polytope such that the length of
every edge of $f$ is equal to $a$ and the length of every
diagonal of $f$ is equal to $b$. Then $f$ is isometric to a
homothetic copy of $S_n$, the standard $n$-dimensional 
cross-polytope in $\Bbb E^n$. In particular, $b=\sqrt{2}a$.
}

\textbf{Proof:} 
As soon as we know distances between every two vertices
of $f$, we may treat $f$ a simplex with $2n$ vertices 
$$
f({\bf e}_1), f(-{\bf e}_1),
f({\bf e}_2), f(-{\bf e}_2), \dots, 
f({\bf e}_n), f(-{\bf e}_n).\eqno(1)
$$ 
In general, a simplex with $2n$ vertices and prescribed
edge lengths is located in $\Bbb E^{2n-1}$.
According to assumptions of Theorem 1, $f$ is located in 
$\Bbb E^{2n-2}$ and,
thus, its $(2n-1)$-dimensional volume is equal to 0.

Let's make use of the Cayley--Menger formula for the 
$k$-dimensional volume, $\mbox{vol}_k$, of a simplex with 
$k+1$ vertices in $\Bbb E^k$ (see, e.g., \cite[p. 98]{Bl53}):
$$
(-1)^{k+1}2^k(k!)^2(\mbox{vol}_k)^2=
\left|
\begin{array}{cccccc}
0     & 1         & 1        & 1        & \dots & 1          \\
1     & 0         & d_{12}^2 & d_{13}^2 & \dots & d_{1,k+1}^2 \\
1     & d_{21}^2  & 0        & d_{23}^2 & \dots & d_{2,k+1}^2 \\
1     & d_{31}^2  & d_{32}^2 & 0        & \dots & d_{3,k+1}^2 \\
\dots & \dots     & \dots    & \dots    & \dots & \dots      \\
1     & d_{k+1,1}^2  & d_{k+1,2}^2 & d_{k+1,3}^2 & \dots & 0          
\end{array}
\right|.
\eqno(2)
$$
Here we assume that the vertices are labeled with
numbers from 1 to $k+1$ and $d_{ij}$ is the Euclidean
distance between $i$-th and $j$-th vetrices.

Enumerating the vertices of $f$ according to (1)
and taking into account that $\mbox{vol}_{2n-1}=0$, 
we obtain from (2)
$$\left|
\begin{array}{ccccc}
0     & 
\begin{array}{cc}
        1 & 1
        \end{array}        & 
\begin{array}{cc}
        1 & 1
        \end{array}        &
        \dots              & 
\begin{array}{cc}
        1 & 1
        \end{array} \\
\begin{array}{cc}
        1 \\ 1
        \end{array}        & 
\fbox{$\begin{array}{cc}
        0 & b^2 \\
        b^2 & 0
        \end{array}$}        &
\begin{array}{cc}
        a^2 & a^2 \\
        a^2 & a^2
        \end{array}        &
\begin{array}{c}
       \dots \\
       \dots 
        \end{array}        &
\begin{array}{cc}
        a^2 & a^2 \\
        a^2 & a^2
        \end{array} \\
\begin{array}{cc}
        1 \\ 1
        \end{array}        & 
\begin{array}{cc}
        a^2 & a^2 \\
        a^2 & a^2
        \end{array}        &
\fbox{$\begin{array}{cc}
        0 & b^2 \\
        b^2 & 0
        \end{array}$}        &
\begin{array}{c}
       \dots \\
       \dots 
        \end{array}        &
\begin{array}{cc}
        a^2 & a^2 \\
        a^2 & a^2
        \end{array} \\
\cdots & \cdots\cdots & \cdots\cdots & \ddots & \cdots\cdots \\
\begin{array}{cc}
        1 \\ 1
        \end{array}        & 
\begin{array}{cc}
        a^2 & a^2 \\
        a^2 & a^2
        \end{array}        &
\begin{array}{cc}
        a^2 & a^2 \\
        a^2 & a^2
        \end{array}        &
\begin{array}{c}
       \dots \\
       \dots 
        \end{array}        &        
\fbox{$\begin{array}{cc}
        0 & b^2 \\
        b^2 & 0
        \end{array}$}        
\end{array}
\right|=0.
\eqno(3)
$$
Note that the matrix in (3) contains $n$ blocks of the form
$$
\left[
\begin{array}{cc}
0   & b^2 \\
b^2 & 0
\end{array}
\right].
$$
In (3), these blocks are boxed for better visibility.

Denote the determinant in (3) by $D_n$.
We are going to compute $D_n$ and, thus, replace (3) 
by an explicit relation involving $a$ and $b$.

Multiply the first row of $D_n$ by $(-a^2)$ and add the result
to every other row. This yields
$$
D_n=\left|
\begin{array}{ccccc}
0     & 
\begin{array}{cc}
        1\hphantom{b^2} & \hphantom{b^2}1
        \end{array}        & 
\begin{array}{cc}
        1\hphantom{b^2} & \hphantom{b^2}1
        \end{array}        &
        \dots              & 
\begin{array}{cc}
        1\hphantom{b^2} & \hphantom{b^2}1
        \end{array} \\
\begin{array}{cc}
        1 \\ 1
        \end{array}        & 
\fbox{$\begin{array}{cc}
        -a^2 & b^2-a^2 \\
        b^2-a^2 & -a^2
        \end{array}$}        &
\begin{array}{cc}
        0\hphantom{b^2} & \hphantom{b^2}0 \\
        0\hphantom{b^2} & \hphantom{b^2}0
        \end{array}        &
\begin{array}{c}
       \dots \\
       \dots 
        \end{array}        &
\begin{array}{cc}
        0\hphantom{b^2} & \hphantom{b^2}0 \\
        0\hphantom{b^2} & \hphantom{b^2}0
        \end{array} \\
\begin{array}{cc}
        1 \\ 1
        \end{array}        & 
\begin{array}{cc}
        0\hphantom{b^2} & \hphantom{b^2}0 \\
        0\hphantom{b^2} & \hphantom{b^2}0
        \end{array}        &
\fbox{$\begin{array}{cc}
        -a^2 & b^2-a^2 \\
        b^2-a^2 & -a^2
        \end{array}$}        &
\begin{array}{c}
       \dots \\
       \dots 
        \end{array}        &
\begin{array}{cc}
        0\hphantom{b^2} & \hphantom{b^2}0 \\
        0\hphantom{b^2} & \hphantom{b^2}0
        \end{array} \\
\cdots & \cdots\cdots & \cdots\cdots & \ddots & \cdots\cdots \\
\begin{array}{cc}
        1 \\ 1
        \end{array}        & 
\begin{array}{cc}
        0\hphantom{b^2} & \hphantom{b^2}0 \\
        0\hphantom{b^2} & \hphantom{b^2}0
        \end{array}        &
\begin{array}{cc}
        0\hphantom{b^2} & \hphantom{b^2}0 \\
        0\hphantom{b^2} & \hphantom{b^2}0
        \end{array}        &
\begin{array}{c}
       \dots \\
       \dots 
        \end{array}        &        
\fbox{$\begin{array}{cc}
        -a^2 & b^2-a^2 \\
        b^2-a^2 & -a^2
        \end{array}$}        
\end{array}
\right|.\eqno(4)
$$
In (4), subtract the second row from the third,
the forth from the fifth, \dots, and the $(2n)$th 
from the $(2n+1)$th. We obtain
$$
D_n=\left|
\begin{array}{ccccc}
0     & 
\begin{array}{cc}
        1\hphantom{b^2} & \hphantom{b^2}1
        \end{array}        & 
\begin{array}{cc}
        1\hphantom{b^2} & \hphantom{b^2}1
        \end{array}        &
        \dots              & 
\begin{array}{cc}
        1\hphantom{b^2} & \hphantom{b^2}1
        \end{array} \\
\begin{array}{cc}
        1 \\ 1
        \end{array}        & 
\fbox{$\begin{array}{cc}
        -a^2 & b^2-a^2 \\
       \hphantom{1} b^2 & -b^2
        \end{array}$}        &
\begin{array}{cc}
        0\hphantom{b^2} & \hphantom{b^2}0 \\
        0\hphantom{b^2} & \hphantom{b^2}0
        \end{array}        &
\begin{array}{c}
       \dots \\
       \dots 
        \end{array}        &
\begin{array}{cc}
        0\hphantom{b^2} & \hphantom{b^2}0 \\
        0\hphantom{b^2} & \hphantom{b^2}0
        \end{array} \\
\begin{array}{cc}
        1 \\ 1
        \end{array}        & 
\begin{array}{cc}
        0\hphantom{b^2} & \hphantom{b^2}0 \\
        0\hphantom{b^2} & \hphantom{b^2}0
        \end{array}        &
\fbox{$\begin{array}{cc}
        -a^2 & b^2-a^2 \\
        \hphantom{1}b^2 & -b^2
        \end{array}$}        &
\begin{array}{c}
       \dots \\
       \dots 
        \end{array}        &
\begin{array}{cc}
        0\hphantom{b^2} & \hphantom{b^2}0 \\
        0\hphantom{b^2} & \hphantom{b^2}0
        \end{array} \\
\cdots & \cdots\cdots & \cdots\cdots & \ddots & \cdots\cdots \\
\begin{array}{cc}
        1 \\ 1
        \end{array}        & 
\begin{array}{cc}
        0\hphantom{b^2} & \hphantom{b^2}0 \\
        0\hphantom{b^2} & \hphantom{b^2}0
        \end{array}        &
\begin{array}{cc}
        0\hphantom{b^2} & \hphantom{b^2}0 \\
        0\hphantom{b^2} & \hphantom{b^2}0
        \end{array}        &
\begin{array}{c}
       \dots \\
       \dots 
        \end{array}        &        
\fbox{$\begin{array}{cc}
        -a^2 & b^2-a^2 \\
        \hphantom{1}b^2 & -b^2
        \end{array}$}        
\end{array}
\right|.\eqno(5)
$$
In (5), add the second column to the third,
the forth to the fifth, \dots, and the $(2n)$th 
to the $(2n+1)$th. We get
$$
D_n=\left|
\begin{array}{ccccc}
0     & 
\begin{array}{cc}
        1\hphantom{b^2} & \hphantom{b^2}2
        \end{array}        & 
\begin{array}{cc}
        1\hphantom{b^2} & \hphantom{b^2}2
        \end{array}        &
        \dots              & 
\begin{array}{cc}
        1\hphantom{b^2} & \hphantom{b^2}2
        \end{array} \\
\begin{array}{cc}
        1 \\ 0
        \end{array}        & 
\fbox{$\begin{array}{cc}
        -a^2 & b^2-2a^2 \\
       \hphantom{1} b^2 & 0
        \end{array}$}        &
\begin{array}{cc}
        0\hphantom{b^2} & \hphantom{b^2}0 \\
        0\hphantom{b^2} & \hphantom{b^2}0
        \end{array}        &
\begin{array}{c}
       \dots \\
       \dots 
        \end{array}        &
\begin{array}{cc}
        0\hphantom{b^2} & \hphantom{b^2}0 \\
        0\hphantom{b^2} & \hphantom{b^2}0
        \end{array} \\
\begin{array}{cc}
        1 \\ 0
        \end{array}        & 
\begin{array}{cc}
        0\hphantom{b^2} & \hphantom{b^2}0 \\
        0\hphantom{b^2} & \hphantom{b^2}0
        \end{array}        &
\fbox{$\begin{array}{cc}
        -a^2 & b^2-2a^2 \\
        \hphantom{1}b^2 & 0
        \end{array}$}        &
\begin{array}{c}
       \dots \\
       \dots 
        \end{array}        &
\begin{array}{cc}
        0\hphantom{b^2} & \hphantom{b^2}0 \\
        0\hphantom{b^2} & \hphantom{b^2}0
        \end{array} \\
\cdots & \cdots\cdots & \cdots\cdots & \ddots & \cdots\cdots \\
\begin{array}{cc}
        1 \\ 0
        \end{array}        & 
\begin{array}{cc}
        0\hphantom{b^2} & \hphantom{b^2}0 \\
        0\hphantom{b^2} & \hphantom{b^2}0
        \end{array}        &
\begin{array}{cc}
        0\hphantom{b^2} & \hphantom{b^2}0 \\
        0\hphantom{b^2} & \hphantom{b^2}0
        \end{array}        &
\begin{array}{c}
       \dots \\
       \dots 
        \end{array}        &        
\fbox{$\begin{array}{cc}
        -a^2 & b^2-2a^2 \\
        \hphantom{1}b^2 & 0
        \end{array}$}        
\end{array}
\right|.\eqno(6)
$$
Expand the determinant in (6) along the third row,
the fifth row, \dots, and the $(2n+1)$th row.
The result is
$$
D_n=(-1)^n b^{2n}\left|
\begin{array}{ccccc}
0 & 2 & 2 & \dots & 2 \\
1 & \fbox{$b^2-2a^2$} & 0 & \dots & 0 \\
1 & 0 & \fbox{$b^2-2a^2$} & \dots & 0 \\
\cdots & \cdots & \cdots & \ddots & \cdots \\
1 & 0 & 0 & \dots & \fbox{$b^2-2a^2$} 
\end{array}
\right|. \eqno(7)
$$
Note that the size of the matrix in (7) is reduced to $n+1$.

If $b^2-2a^2=0$, the determinant in (7) is equal to zero
(e.g., because the second row is equal to the third).
If $b^2-2a^2\neq 0$, rewrite (7) in the form
$$
D_n=2\frac{(-1)^n b^{2n}}{b^2-2a^2}\left|
\begin{array}{ccccc}
0 & 1 & 1 & \dots & 1 \\
b^2-2a^2 & \fbox{$b^2-2a^2$} & 0 & \dots & 0 \\
b^2-2a^2 & 0 & \fbox{$b^2-2a^2$} & \dots & 0 \\
\cdots & \cdots & \cdots & \ddots & \cdots \\
b^2-2a^2 & 0 & 0 & \dots & \fbox{$b^2-2a^2$} 
\end{array}
\right|.\eqno(8)
$$
In (8), subtract the second, the third, \dots, and the $(n+1)$st 
column from the first column.
This yields
$$
D_n=2\frac{(-1)^n b^{2n}}{b^2-2a^2}\left|
\begin{array}{ccccc}
-n & 1 & 1 & \dots & 1 \\
0 & \fbox{$b^2-2a^2$} & 0 & \dots & 0 \\
0 & 0 & \fbox{$b^2-2a^2$} & \dots & 0 \\
\cdots & \cdots & \cdots & \ddots & \cdots \\
0 & 0 & 0 & \dots & \fbox{$b^2-2a^2$} 
\end{array}
\right|= 2nb^{2n}(2a^2-b^2)^{n-1}.
$$
Hence, (3) is equivalent to $b^2-2a^2=0$. 
This means that the cross-polytope $f$ is congruent to a homothetic copy
of the standard $n$-dimensional cross-polytope $S_n$ in $\Bbb E^n$. Q.E.D.

\textbf{2. Two-distance preserving mappings.}
A mapping $g:\Bbb E^n \to \Bbb E^m$ is said to be \textit{unit 
distance preserving} if, for all $x,y\in \Bbb E^n$,
the equality $|x-y|=1$ implies the equality $|g(x)-g(y)|=1$.
Here by $|x|$ we denote the Euclidean norm of a vector 
$x\in\Bbb E^n$.

In 1953, F.S. Beckman and D.A. Quarles  \cite{BQ53}
proved that, \textit{for $n\geqslant 2$, every unit distance preserving 
mapping $g:\Bbb E^n \to \Bbb E^n$ is an isometry of $\Bbb E^n$}
(i.e., $g$ preserves all distances).
Since that time, the problem `does a unit distance preserving 
mapping is an isometry' was studied for spaces of various types
(hyperbolic \cite{Ku80}, Banach \cite{Ra07}, 
$\Bbb Q^n$ \cite{Za05}, just  to name a few). 

In 1985, B.V. Dekster \cite{De85}
found a mapping $g:\Bbb E^2 \to \Bbb E^6$, which is 
unit distance preserving but is not an isometry.
That example motivated geometers to look for other
conditions that make the statement `every unit distance 
preserving mapping $g:\Bbb E^n \to \Bbb E^m$ is an isometry' 
correct even if $m\neq n$.
One of the possible sets of such conditions uses the notions 
of cable and strut defined as follows.

Given a mapping $g:\Bbb E^n \to \Bbb E^m$, a positive
real number $c$ is called a \textit{cable} of $g$ if, for all
$x,y\in\Bbb E^n$, the equality $|x-y|=c$ implies
$|g(x)-g(y)|\leqslant c$ and a positive real number
$s$ is called a \textit{strut} of $g$ if, for all
$x,y\in\Bbb E^n$, the equality $|x-y|=s$ implies
$|g(x)-g(y)|\geqslant s$.

In 1999, K. Bezdek and R. Connelly \cite{BC99} proved that
\textit{if $n\geqslant 2$, $c$ is a cable of a mapping 
$g:\Bbb E^n \to \Bbb E^m$, $s$ is a strut of $g$ and 
$c/s<(\sqrt{5}-1)/2$, then $g$ is an isometry.}

As a corollary of Theorem 1, we prove the following theorem.

\textbf{Theorem 2.} \textit{Let $n\geqslant 6$ and 
$0\leqslant m\leqslant 2n-2$ be integers,
let $A$ and $B$ be positive real numbers, and let 
$g:\Bbb E^n \to \Bbb E^m$ be a mapping such that, for all $x,y\in\Bbb E^n$,
the equality $|x-y|=A$ implies $|g(x)-g(y)|=A$ and
the equality $|x-y|=\sqrt{2}A$ implies $|g(x)-g(y)|=B$.
Then $g$ is an isometry.}

\textbf{Proof:} 
First, observe that 
\textit{$B$ is necessarily equal to $\sqrt{2}A$}.

In fact, given $x,y\in\Bbb E^n$ such that
$|x-y|=\sqrt{2}A$, find an $n$-dimensional 
cross-polytope $P$ in $\Bbb E^n$ with the following properties:

(i) $P$ is congruent to a homothetic copy of the standard 
$n$-dimensional cross-polytope $S_n$ in $\Bbb E^n$ with 
the scale factor $A/\sqrt{2}$;

(ii) $x$ and $y$ belong to the vertex set of $P$;

(iii)  the straight line segment with the endpoints $x$ 
and $y$ is a diagonal of $P$.

Since $g:\Bbb E^n \to \Bbb E^m$ and $m\leqslant 2n-2$,
Theorem 1 yields that the image of $P$ under the mapping $g$
is congruent to $P$. In particular, this means that
$|g(x)-g(y)|=\sqrt{2}A$ and, thus, $B=\sqrt{2}A$.

Now, let's prove that \textit{the real number $c=2A/\sqrt{n}$ 
is a cable of the mapping $g$,} i.e., let's prove that
if $v_1,v_2\in\Bbb E^n$ are such that 
$|v_1-v_2|=c=2A/\sqrt{n}$ 
then $|g(v_1)-g(v_2)|\leqslant c$.

Let $L$ be the $(n-1)$-dimensional plane in $\Bbb E^n$ that
passes thought the point $(v_1+v_2)/2$ and is orthogonal 
to the vector $v_1-v_2$.
Let $v_3, v_4, \dots , v_{n+2}$ be the vertices of an
$(n-1)$-dimensional regular simplex in $L$ with edge lengths
$\sqrt{2}A$ and circumcenter at the point $(v_1+v_2)/2$.
The latter means that the point $(v_1+v_2)/2$ is the center
of an $(n-2)$-dimensional sphere which passes through all 
the vertices $v_3, v_4, \dots , v_{n+2}$.
Since the radius of this $(n-2)$-sphere is equal to 
$R=A\sqrt{1-1/n}$ (see, e.g., \cite[pp. 294--295]{Co48}),
Pythagora's Theorem gives
$|v_i-v_j|=A$ for all $i=1,2$ and $j=3,4,\dots, n+2$.

It follows from conditions of Theorem 2 
and the relation $B=\sqrt{2}A$ that, for every $i=1,2$,
the $n$-simplex $\Delta_i$
with vertices $v_i, v_3, v_4, \dots , v_{n+2}$ 
is congruent to the $n$-simplex with vertices
$g(v_i), g(v_3), g(v_4), \dots , g(v_{n+2})$. 
For short, denote the latter simplex by $g(\Delta_i)$.
Denote by $\lambda$ the $(n-1)$-dimensional plane 
containing the points $g(v_3), g(v_4), \dots , g(v_{n+2})$. 

Since $g(\Delta_1)$ is congruent to $\Delta_1$, it follows
that $m\geqslant n$. 

If $m=n$, the non-degenerate $n$-simplices $g(\Delta_1)$ 
and $g(\Delta_2)$ lie either in the same half-space of 
$\Bbb E^m$ determined by $\lambda$ 
(in this case $|g(v_1)-g(v_2)|=0<c$) 
or in the different half-spaces of $\Bbb E^m$ determined 
by $\lambda$ (in this case 
$|g(v_1)-g(v_2)|= |v_1-v_2|=c$). In both cases, 
$|g(v_1)-g(v_2)|\leqslant c$. 

If $m>n$, the $n$-simplices $g(\Delta_1)$ and $g(\Delta_2)$
may be obtained from $\Delta_1$ and $\Delta_2$ in two steps:
first, we apply to $\Delta_1$ and $\Delta_2$ such an isometry
$h:\Bbb E^n\to\Bbb E^m$ that $h(v_j)=g(v_j)$ for all 
$j=1,3,4,\dots, n+2$ and then rotate the simplex $h(\Delta_2)$ 
around $\lambda$ in such a way that $h(v_2)$ coincides with $g(v_2)$.
From this description we conclude that 
$|g(v_1)-g(v_2)|\leqslant |h(v_1)-h(v_2)|=|v_1-v_2|=c$.
Hence, $c=2A/\sqrt{n}$ is a cable.

Obviously, we may consider $s=\sqrt{2}A$ as a strut of $g$.

Since $n\geqslant 6$, $c/s=\sqrt{2/n}<(\sqrt{5}-1)/2$.
Thus, according to the above-cited theorem by K. Bezdek and R. Connelly,
$g$ is an isometry. Q.E.D.

\bigskip

\noindent{Victor Alexandrov}

\noindent\textit{Sobolev Institute of Mathematics}

\noindent\textit{Koptyug ave., 4}

\noindent\textit{Novosibirsk, 630090, Russia}

and

\noindent\textit{Department of Physics}

\noindent\textit{Novosibirsk State University}

\noindent\textit{Pirogov str., 2}

\noindent\textit{Novosibirsk, 630090, Russia}

\noindent\textit{e-mail: alex@math.nsc.ru}

\bigskip

\noindent{Submitted: April 17, 2015}

\end{document}